%% file: agt-3-37.tex
\documentclass{gtart}
\input agtout

\lognumber{37}
\volumenumber{3}
\volumeyear{2003}
\papernumber{37}
\published{22 October 2003}
\pagenumbers{1079}{1088}
\received{17 March 2002}
\revised{25 August 2003}
\accepted{24 September 2003}

\usepackage{amsmath,amssymb}
\def\qedhere{\eqno{\qed}}

\newtheorem{thm}{Theorem}[section]
\newtheorem{cor}[thm]{Corollary}

\theoremstyle{remark}
\newtheorem{remark}[thm]{Remark}

\theoremstyle{remark}

\newcommand{\gp}[1]{\left\langle\,#1\,\right\rangle}
\newcommand{\id}{\mathit{id}}
\newcommand{\C}{\mathcal{C}}

\begin{document}
\title{Rigidity of graph products of groups}
\author{David G. Radcliffe}
\address{1924 Ford Parkway \#10\\Saint Paul MN 55116, USA}
\email{radcliffe@xirr.com}

\begin{abstract}
We show that if a group can be represented as a graph product of finite directly indecomposable groups, then this
representation is unique.
\end{abstract}

 \maketitle

\section{Introduction}

Given a simple graph with nontrivial groups as vertices, a group is formed by taking the free product of the vertex
groups, with added relations implying that elements of adjacent groups commute. This group is said to be the
\emph{graph product} of the vertex groups. If the graph is discrete then the graph product is the free product of the
vertex groups; while if the graph is complete then the graph product is the restricted direct product\footnote{By
``restricted'' we mean that all but finitely many entries are the identity.} of the vertex groups. Graph products were
first defined in Elisabeth Green's Ph.D. thesis \cite{elisabethgreen}, and have been studied by other authors
\cite{MR96a:20052, MR2000k:20056, MR2002k:20074}.

Important special cases of graph products arise when we specify the vertex groups. If all vertex groups are infinite
cyclic, then the graph product is called a \emph{graph group} or a \emph{right-angled Artin group}. Graph groups have
been studied by many authors \cite{MR88e:20033, MR90m:20043, MR95g:20041}. If all vertex groups have order two, then the
graph product is called a \emph{right-angled Coxeter group}. These groups were first studied by Ian Chiswell
\cite{MR89g:20053}, and they have been studied by many other authors \cite{MR2001m:20056, MR2002e:58039, MR2000d:20069}.

In this article we investigate the question of uniqueness for graph product decompositions. Carl Droms \cite{MR88e:20033}
proved that two graph products of infinite cyclic groups are isomorphic if and only if their graphs are isomorphic.
Elisabeth Green \cite{elisabethgreen} proved that if a group can be represented as a graph product of cyclic groups of prime
order, then this representation is unique. This result was extended to primary cyclic groups by the present author
\cite{davidradcliffe}. Our main result is the following: If a group can be represented as a graph product of directly
indecomposable finite groups, then this representation is unique.

\section{Graphs and modular partitions}

A graph is an ordered pair of sets $(V,E)$ where $E$ is a set of two-element subsets of $V$. Elements of $V$ are called
\emph{vertices}, and elements of $E$ are called \emph{edges}. For the remainder of this paper we shall assume that $V$
is finite. A clique is a maximal complete subgraph, or (by abuse of terminology) the set of vertices of a maximal
complete subgraph.

A \emph{module} of a graph $(V,E)$ is a subset $X$ of $V$ such that for every $v \in V-X$, either $v$ is adjacent to
every element of $X$ or $v$ is adjacent to no element of $X$. A \emph{modular partition} is a partition of $V$ into
non-empty modules. A modular partition induces a quotient graph $(\bar{V},\bar{E})$ where $\bar{V}$ is the set of
partition classes and $\{\alpha,\beta\}\in\bar{E}$ if and only $\{u,v\}\in E$ for some (and hence for all) $u\in\alpha$
and $v\in\beta$. We may regard $(\bar{V},\bar{E})$ as a compressed version of the original graph. Given the quotient
graph and the subgraphs induced by the partition classes, it is possible to reconstruct the original graph. For this
reason, modular partitions have been studied extensively by computer scientists \cite{MR2000d:05117}.

We say that a graph $(V,E)$ is \emph{$T_0$} if no edge is a module. This means that for all $\{u,v\} \in E$, there
exists $w \in V-\{u,v\}$ so that $w$ is adjacent to $u$ or $v$ but not both. A graph is $T_0$ if and only if vertices
are distinguished by the cliques to which they belong. That is, a graph is $T_0$ if and only if the following condition
holds: for every pair of distinct vertices, there exists a clique which contains exactly one of them.

Let us say that two vertices are \emph{equivalent} if they cannot be distinguished by the cliques. Then the set of
equivalence classes is a modular partition. The quotient graph resulting from this partition satisfies the $T_0$
condition, and it will be called the \emph{$T_0$ quotient}.

Similarly, a graph $(V,E)$ is \emph{$T_1$} if for all $\{u,v\} \in E$ there exists $w \in V-\{u,v\}$ so that $\{u,w\}
\in E$ and $\{v,w\} \not\in E$. Equivalently, a graph is $T_1$ if and only if every vertex is the intersection of the
set of cliques to which it belongs. Note that this condition is stronger than the $T_0$ condition.

\section{Graph products of groups}

Let $\Gamma = (V,E)$ be a graph, and let $\{G_v\}_{v\in V}$ be a collection of groups which is indexed by the vertex
set of $\Gamma$. We say that $(\Gamma, G_v)$ is a \emph{graph of groups.}\footnote{This differs from the usual
definition, which has vertex groups and edge groups, together with monomorphisms from the edge groups to the vertex
groups \cite{MR94j:20034}.} Two graphs of groups, $(\Gamma, G_v)$ and $(\Gamma', G'_v)$, are \emph{isomorphic} if
there exists a graph isomorphism $\phi\colon\Gamma\rightarrow\Gamma'$ so that $G_v$ and $G'_{\phi(v)}$ are isomorphic
for all $v\in V$.

The \emph{graph product} $G$ of a graph of groups is the quotient of the free product of the vertex groups by the
normal subgroup generated by all commutators of elements taken from pairs of adjacent groups. That is, $G=F/N$ where $F
= \coprod_{v\in V} G_v$ and $N$ is the normal closure in $F$ of $$\{g^{-1}h^{-1}gh \ \colon g \in G_u, h \in G_v,
\{u,v\}\in E\}.$$ The canonical monomorphism from $G_v$ to $F$ induces a monomorphism from $G_v$ to $G$. We may thus
identify each vertex group $G_v$ with its image in $G$, in which case we say that $G$ is an \emph{internal} graph
product.

The graph product can also be described in terms of generators and relations. Choose a presentation $(\gamma_v;
\rho_v)$ for each vertex group, so that the generating sets $\gamma_v$ are pairwise disjoint. Then $G$ has a
presentation $(\bigcup \gamma_v ; \bigcup\rho_v \cup \sigma)$ where $\sigma = \{a^{-1}b^{-1}ab\colon a\in \gamma_u, b
\in \gamma_v, \{u,v\}\in E\}$.

If $A$ is a subset of $V$, then we denote by $\Gamma_A$ the subgraph of $\Gamma$ that is induced by $A$, and we denote
by $G(A)$ the subgroup of $G$ that is generated by $\bigcup_{a\in A} G_a$.

\begin{thm} \cite{elisabethgreen, davidradcliffe} If $A$ is a subset of $V$ then $G(A)$ is the internal graph product of $(\Gamma_A, G_a)$.
\qed \end{thm}

\begin{cor}
\label{sum} If $A$ is complete then $G(A)$ is the (restricted) direct product of the $G_a$, and if $A$ is discrete then
$G(A)$ is the free product of the $G_a$. \qed
\end{cor}

\begin{thm}
If $A \subseteq V$ then there is a homomorphism $\rho_A \colon G \rightarrow G(A)$ so that $\rho_A(x)=x$ for all $x\in
G(A)$ and $\rho_A(x)=1$ for all $x\in G(V-A)$. We call $\rho_A$ a \emph{retraction homomorphism}.\label{retraction}
\end{thm}
\begin{proof}
For each $a\in A$ let $h_a \colon G_a \rightarrow G(A)$ be the inclusion homomorphism, and for $b \in V-A$ let $h_b
\colon G_b \rightarrow G(A)$ be the trivial homomorphism. Then there exists a homomorphism $\rho_A$ which extends $h_v$
for all $v \in V$. It is clear that $\rho_A(x)=x$ for all $x\in G(A)$, since $\rho_A(x)=x$ for all $x \in \bigcup_{a\in
A} G_a$, and likewise that $\rho_A(x)=1$ for all $x\in G(V-A)$.
\end{proof}

\begin{thm} \label{cup}
If $A$ and $B$ are subsets of $V$ then $G(A\cup B) = \gp{G(A) \cup G(B)}$.
\end{thm}
\proof
$G(A) = \gp{\bigcup_{a\in A} G_a}$ and $G(B) = \gp{\bigcup_{b \in B} G_b}$, thus $$\gp{G(A) \cup G(B)} =
\gp{\bigcup_{c\in A\cup B} G_c} = G(A\cup B).\qedhere$$

\begin{thm} \label{cap}
If $A$ and $B$ are subsets of $V$ then $G(A) \cap G(B) = G(A\cap B)$.
\end{thm}

\begin{proof}
It is clear that $G(A\cap B) \subseteq G(A) \cap G(B)$. For the reverse inclusion, let $\rho_A \colon G \rightarrow
G(A)$ be the retraction homomorphism of Theorem \ref{retraction}, and let $x \in G(A) \cap G(B)$.  It remains to prove
that $x \in G(A \cap B)$. If $b \in A \cap B$ then $\rho_A(y)=y$ for all $y\in G_b$. If $b \in B - A$ then $\rho_A(y) =
1$ for all $x\in G_b$. In either case $\rho_A(y) \in G(A\cap B)$ for all $y\in G_b$ and all $b\in B$. Therefore
$\rho_A(x) \in G(A\cap B)$. But $\rho_A(x) = x$ since $x\in G(A)$. Therefore $x\in G(A\cap B)$ as claimed.
\end{proof}

\begin{thm}
Let $A,B\subseteq V$ be complete. If $x\in G(A)$ and $x$ is conjugate to an element $y\in G(B)$, then $x\in G(A\cap B)$.
\label{conj}
\end{thm}

\begin{proof}
Let $\rho \colon G \rightarrow G(V-B)$ be the retraction homomorphism of Theorem \ref{retraction}. Then $\rho(y)=1$, so
$\rho(x)=1$ as well, since the kernel is a normal subgroup.

By Corollary \ref{sum}, we may express $x$ uniquely as $x = \prod_{a \in A} x_a$, where $x_a \in G_a$ for all $a \in
A$. Then $\rho(x) = \prod_{a \in A} \rho(x_a) = \prod_{a \in A-B} x_a$. But $\rho(x)=1$, so $x_a = 1$ for all $a \in
A-B$. Therefore $x \in G(A\cap B)$.
\end{proof}

\begin{cor} \label{conj2}
If $A,B\subseteq V$ are complete and $G(A)$ is conjugate to $G(B)$ then $A=B$. \qed
\end{cor}

The proof of the following theorem is left to the reader.

\begin{thm}
Let $\bar{\Gamma} = (\bar{V},\bar{E})$ be the quotient graph resulting from a modular partition of $\Gamma = (V,E)$.
Then $G$ is the graph product of $(\bar{\Gamma}, G(A))$, where $A$ varies over the modules of $\Gamma$.  \label{t0q}
\qed
\end{thm}

Since the partition of a graph into its components (or co-components) is modular, we obtain the following corollary.
(Recall that a co-component of a graph is a component of the complement.)

\begin{cor}
If the components of $\Gamma$ are $A_1, \ldots, A_n$ then $G \cong \coprod_i G(A_i)$. If the co-components of $\Gamma$
are $B_1, \ldots, B_m$ then $G \cong \bigoplus_i G(B_i)$. \qed
\end{cor}

We also require the following result, which is proved in \cite{elisabethgreen}.

\begin{thm}
For every finite subgroup $F$ of $G$ there exists a complete subgraph $C$ so that $F$ is conjugate to a subgroup of
$G(C)$. \qed \label{finsub}
\end{thm}

\begin{cor}
If all vertex groups are finite, then $F$ is a maximal finite subgroup of $G$ if and only if there exists a clique $C$
so that $F$ is conjugate to $G(C)$. \label{maxfin}
\end{cor}

\begin{proof}
Let $F$ be a maximal finite subgroup of $G$. By the previous theorem, there exists a clique $C$ so that $F$ is
conjugate to a subgroup of $G(C)$. However, $G(C)$ itself is a finite subgroup of $G$. Therefore $F$ is conjugate to
$G(C)$.

Conversely, let $C$ be a clique, and let $F$ be a conjugate of $G(C)$. Let $F'$ be a finite subgroup of $G$ so that $F'
\supseteq F$. By the previous theorem there exists a clique $D$ so that $F'$ is conjugate to a subgroup of $G(D)$.
Therefore $G(C)$ is conjugate to a subgroup of $G(D)$. It follows from Theorem \ref{conj} that $C$ is a subset of $D$.
But $C$ is a clique, hence $C=D$ and $F'=F$. Consequently $F$ is a maximal finite subgroup.
\end{proof}

\begin{remark} An alternate proof of this corollary can be obtained by considering the action of $G$ on the ${\rm CAT}(0)$
cube complex defined by John Meier and other authors \cite{MR97b:20055,MR2000i:20068,MR2002k:20074}.
Any finite group acting on cellularly on a CAT(0) complex fixes some cell. Since stabilizers of cubes in this complex are
conjugates of the groups G(C), the corollary follows.
\end{remark}

\section{Conjugacy classes of finite subgroups}

Let $G$ be the internal graph product of $(\Gamma, G_v)$. We assume for the remainder of this article that $\Gamma =
(V,E)$ is a finite graph, and that each vertex group $G_v$ is finite.

Let $\mathcal{F}$ denote the set of conjugacy classes of finite subgroups of G. We write $[F]$ to denote the set of
subgroups of $G$ which are conjugate to a given finite subgroup $F$. We define a partial ordering on $\mathcal{F}$ as
follows: If $A$ and $B$ are finite subgroups of $G$, then $[A] \preceq [B]$ if and only if there exists $g \in G$ so
that $A \subseteq gBg^{-1}$.

\begin{thm} The relation $\preceq$ is a well-defined partial ordering on $\mathcal{F}$.
\end{thm}

\begin{proof}

\emph{$\preceq$ is well defined:} Let $A$ and $B$ be finite subgroups of $G$, and suppose that there exists $g\in G$ so
that $A \subseteq gBg^{-1}$. Let $A'$ and $B'$ be subgroups of $G$ which are conjugate to $A$ and $B$ respectively.
There exist $h,k \in G$ so that $A' = hAh^{-1}$ and $B' = kBk^{-1}$.  Now $A' \subseteq h(gBg^{-1})h^{-1} =
hgk^{-1}B'kg^{-1}h^{-1}$, so $A' \subseteq mB'm^{-1}$ where $m=hgk^{-1}$.

\emph{$\preceq$ is transitive:} Let $A,B,C$ be finite subgroups of $G$ so that $[A] \preceq [B]$ and $[B] \preceq [C]$.
There exist $g,h \in G$ so that $A \subseteq gBg^{-1}$ and $B \subseteq hCh^{-1}$. Then $A \subseteq ghC(gh)^{-1}$,
hence $[A] \preceq [C]$.

\emph{$\preceq$ is irreflexive:} Let $A$ and $B$ be finite subgroups of $G$ so that $[A] \preceq [B]$ and $[B] \preceq
[A]$. Then there exist $g,h$ so that $A \subseteq gBg^{-1}$ and $B \subseteq hAh^{-1}$. Since $A$ and $B$ are finite,
it follows that $|A| = |B|$ and $A=gBg^{-1}$. Therefore $[A] = [B]$.
\end{proof}

Recall that if $A$ is a subset of a partially ordered set $(X,\leq)$ then the \emph{least upper bound} of $A$, denoted
$\bigvee A$, is an element $x \in X$ so that $a \leq x$ for all $a \in A$, and if $a \leq y$ for all $a \in A$ then $x
\leq y$. The least upper bound is unique when it exists. Similarly, the \emph{greatest lower bound} of $A$, denoted
$\bigwedge A$, is an element $x \in X$ so that $x \leq a$ for all $a \in A$, and if $y \leq a$ for all $a \in A$ then
$y \leq x$.

\begin{thm}
If $A,B \in \C$ then $[G(A\cap B)] = [G(A)] \wedge [G(B)]$. \label{wedge}
\end{thm}

\begin{proof}
It is obvious that $[G(A\cap B)] \preceq [G(A)]$ and $[G(A\cap B)] \preceq [G(B)]$. Suppose that $F$ is a finite
subgroup of $G$ so that $[F] \preceq [G(A)]$ and $[F] \preceq [G(B)]$. We need to show that $[F] \preceq [G(A\cap B)]$.

We may assume without loss of generality that $F \subseteq G(A)$. If $x \in F$ then $x$ is conjugate to an element of
$G(B)$, hence $x \in G(B)$ by Theorem \ref{conj}. Therefore $F \subseteq G(A) \cap G(B) = G(A\cap B)$ by Theorem
\ref{cap}, so we are done.
\end{proof}

\begin{thm}
Let $A,B \in \C$. If $A \cup B \in \C$ then $[G(A\cup B)] = [G(A)] \vee [G(B)]$. If $A \cup B \not\in \C$, then
$[G(A)]$ and $[G(B)]$ do not have a common upper bound. \label{vee}
\end{thm}

\begin{proof}
Suppose that $A \cup B \in \C$. Then $[G(A\cup B)]$ is an upper bound for $[G(A)]$ and $[G(B)]$. We wish to show that
it is the least upper bound.

Let $[F]$ be another upper bound of $[G(A)]$ and $[G(B)]$. By Theorem \ref{finsub}, there exists $C \in \C$ and $h\in
G$ so that $hFh^{-1} \subseteq G(C)$. Then $C \supseteq A\cup B$ by Theorem \ref{conj}.

Now $hFh^{-1}$ contains conjugates of $G(A)$ and $G(B)$, so $hFh^{-1}$ contains both $G(A)$ and $G(B)$ by Theorem
\ref{conj}. Therefore $hFh^{-1} \supseteq G(A\cup B)$, and hence $[G(A\cup B)] \preceq [F]$.

Now suppose that $A \cup B \not\in \C$. If $[F]$ is an upper bound for $[G(A)]$ and $[G(B)]$, then (since $F$ is
finite) there exists a complete subgraph $D$ so that $F \subseteq G(D)$. By Theorem \ref{conj}, $D$ contains $A \cup
B$. But $A \cup B$ is not complete, and this is a contradiction.
\end{proof}

\section{Uniqueness of graph product decompositions}

Let $\Gamma = (V,E)$ and $\Gamma' = (V',E')$ be finite graphs so that $V \cap V' = \emptyset$. Let $G$ be a group, and
suppose that $G_v$ is a nontrivial finite subgroup of $G$ for each $v \in V \cup V'$. Finally, suppose that $G$ is the
internal graph product of both $(\Gamma, G_v)$ and $(\Gamma', G_{v'})$.

\begin{thm} For each clique $C$ of $\Gamma$ there is a unique clique $C'$ of $\Gamma'$ so that
$[G(C)] = [G(C')]$.\label{cliques-correspond}
\end{thm}

\begin{proof}
If $C$ is a clique of $\Gamma$ then $G(C)$ is a maximal finite subgroup of $G$ by Corollary \ref{maxfin}. Again by
Corollary \ref{maxfin}, $G(C)$ is conjugate to $G(C')$ for some clique $C'$ of $\Gamma'$. Therefore $[G(C)] = [G(C')]$.

Uniqueness of $C'$ follows from Corollary \ref{conj2}.
\end{proof}

\begin{thm} If $\Gamma$ and $\Gamma'$ are $T_1$ then there is a graph isomorphism $\phi \colon
\Gamma \rightarrow \Gamma'$ so that $[G_v] = [G_{\phi(v)}]$ for all $v \in V$. In particular, $(\Gamma, G_v)$ and
$(\Gamma', G_{v'})$ are isomorphic graphs of groups.
\end{thm}

\begin{proof}
Let $v \in V$, and let $\{C_1, \ldots, C_n\}$ be the set of all cliques of $\Gamma$ which contain $v$. Then $\{v\} =
\bigcap_{i=1}^n C_i$ since $\Gamma$ is $T_1$.

For each $i$ there exists a clique $C'_i$ of $\Gamma'$ so that $[G(C_i)] = [G(C'_i)]$ by Theorem
\ref{cliques-correspond}. Now $[G_v] = \bigwedge_i [G(C_i)] = \bigwedge_i [G(C'_i)] = [G(C')]$, where $C' = \bigcap_i
C'_i$. In particular, $C'$ is not empty since $[G(C')]$ is non-trivial.

I claim that $C'$ has only one element. Suppose that $C'$ contains two distinct elements $r',s'$. Since $\Gamma'$ is
$T_1$, there is a clique $D'$ of $\Gamma'$ which contains $r'$ but not $s'$, and by Theorem \ref{cliques-correspond}
there is a corresponding clique $D$ of $\Gamma$ so that $[G(D)] = [G(D')]$.

Now $D \cap \bigcap_{i=1}^n C_i = \emptyset$, since $v\not\in D$. But $r' \in D' \cap \bigcap_{i=1}^n C'_i$, which is a
contradiction. Therefore $C'$ has a unique element $v'=\phi(v)$ as claimed, and so $[G_v] = [G_{v'}]$.

Similarly, there is a function $\phi' \colon V' \rightarrow V$ so that $[G_v] = [G_{\phi(v)}]$ for all $v\in V'$. Then
$\phi' \circ \phi = \id_V$ and $\phi \circ \phi' = \id_{V'}$, so $\phi$ is a bijection and $\phi' = \phi^{-1}$.

If $\{u,v\}\in E$ then there is a clique $C$ of $\Gamma$ so that $\{u,v\} \subseteq C$. Then $\{\phi(u),\phi(v)\}
\subseteq C'$, therefore $\{\phi(u),\phi(v)\} \in E'$. Conversely, if $\{\phi(u),\phi(v)\} \in E'$ then there exists a
clique $C'$ of $\Gamma'$ so that $\{\phi(u),\phi(v)\} \subseteq C'$.  Hence $\{u,v\}\subseteq C$ and so $\{u,v\} \in
E$. Therefore $\phi$ is an isomorphism of graphs, as claimed.
\end{proof}

\begin{thm} If $\Gamma$ and $\Gamma'$ are $T_0$ graphs, then $(\Gamma, G_v)$ and $(\Gamma', G_{v'})$ are
isomorphic graphs of groups.
\end{thm}

\begin{proof}
Let $v \in V$. Let $\{C_1, \ldots, C_n\}$ be the set of cliques of $\Gamma$ which contain $v$, and let $\{D_1, \ldots,
D_m\}$ be the set of cliques of $\Gamma$ which do not contain $v$. It follows from the $T_0$ hypothesis that $\{v\} =
\bigcap_i C_i - \bigcup_i D_i$.

For each $C_i$ and each $D_i$ there are cliques $C'_i$ and $D'_i$ of $\Gamma'$ so that $[G(C_i)] = [G(C'_i)]$ and
$[G(D_i)] = [G(D'_i)]$. Observe that $C - \{v\} = \bigcup_i (C \cap D_i)$. Let $C' = \bigcap_i C'_i$ and $D' = \bigcup
(C' \cap D'_i)$.

Now
$$
\begin{array}{l}
[G(C)] = [G(\bigcap\nolimits_i C_i)] = \bigwedge\nolimits_i [G(C_i)] = \bigwedge\nolimits_i [G(C'_i)] =
[G(\bigcap\nolimits_i C'_i)] =  [G(C')]
\end{array}
$$ and
\renewcommand{\arraystretch}{.5}
$$ \begin{array}{l} [G(C-\{v\})] = [G(\bigcup\nolimits_i (C \cap D_i))] = \bigvee\nolimits_i [G(C \cap D_i)] =
\bigvee\nolimits_i ([G(C)] \wedge
[G(D_i)]) \\
\\
= \bigvee\nolimits_i \left([G(C')] \wedge [G(D'_i)]\right) = \bigvee\nolimits_i[G(C'\cap D'_i)]  =
[G(\bigcup\nolimits_i (C' \cap D'_i))] = [G(D')].
\end{array}$$
\renewcommand{\arraystretch}{1}

Choose $h \in G$ so that $G(C) = hG(C')h^{-1}$. Then $hG(C-\{v\})h^{-1}$ is a subgroup of $G(C')$ that is conjugate to
$G(D')$. Theorem \ref{conj} implies that $hG(C-\{v\})h^{-1} = G(D')$. Therefore,
$$G_v \cong G(C)/G(C-\{v\}) \cong G(C')/G(D') \cong G(C'-D').$$

In particular, $C'-D'$ is nonempty. But the $T_0$ hypothesis prevents $C'-D'$ from having more than one element, since
two elements of $C'-D'$ would belong to the same cliques of $\Gamma'$. So $C'-D'$ has a unique element $v' = \phi(v)$.

In a similar manner, we can associate to each $v' \in V'$ a unique element $v = \phi'(v')$ of $V$. Now $\phi' \circ
\phi(v)$ belongs to the same cliques that $v$ does, so $\phi' \circ \phi = \id_V$, and likewise $\phi \circ \phi' =
\id_{V'}$.

Therefore $\phi$ is a bijection from $V$ to $V'$ such that $G_v \cong G_{\phi(v)}$ for all $v\in V$. It remains only to
prove that $\phi$ is a graph isomorphism. Now, if $\{u,v\} \in E$ then there exists a clique $C$ so that $\{u,v\}
\subset C$. Then $\{\phi(u),\phi(v)\} \subset C'$, where $C'$ is a clique of $\Gamma'$ and $[G(C)] = [G(C')]$.
Therefore $\{\phi(u),\phi(v)\} \in E'$. A similar argument shows that if $\{\phi(u),\phi(v)\} \in E'$ then $\{u,v\} \in
E$. So $\phi$ is a graph isomorphism, and we are done.
\end{proof}

\begin{thm}
If $G_v$ is directly indecomposable for all $v\in V \cup V'$ then $(\Gamma, G_v)$ and $(\Gamma', G_{v'})$ are
isomorphic graphs of groups.
\end{thm}

\begin{proof}
Let $\bar{\Gamma} = (\bar{V},\bar{E})$ and $\bar{\Gamma'} = (\bar{V'},\bar{E'})$ be the $T_0$ quotients of $\Gamma =
(V,E)$ and $\Gamma' = (V',E')$ respectively. Then $G$ is the internal graph product of both $(\bar{\Gamma},G(A))$ and
$(\bar{\Gamma'},G(A'))$.

If $A\in\bar{V}$ then $G(A) = \oplus_{a \in A} G_a$, so $G(A)$ is a finite group. Likewise each $G(A')$ is a finite
group. By the previous theorem, there exists a graph isomorphism $\bar{\phi} \colon \bar{V} \rightarrow \bar{V'}$ so
that $G(A) \cong G(\bar{\phi}(A))$ for all $A \in \bar{V}$.

It is well-known that every finite group has a unique factorization as a direct product of directly indecomposable
groups, up to isomorphism and order of factors \cite{MR95m:20001}. Thus, for each $A \in \bar{V}$ there is a bijection
$\phi_A \colon A \rightarrow \bar{\phi}(A)$ so that $G_v \cong G_{\phi_A(v)}$ for all $v\in A$.

Let $\phi \colon V \rightarrow V'$ be the union of the $\phi_A$'s. Then $\phi$ is clearly a graph isomorphism between
$\Gamma$ and $\Gamma'$, and $G_v \cong G_{\phi(v)}$ for all $v\in V$. Therefore $\phi$ is an isomorphism between the
two graphs of groups.
\end{proof}

\Addresses\recd
\end{document}

%% file: agtout.tex
\def\ifplaintex{\expandafter\ifx\csname documentclass\endcsname\relax}

\def\gtp{{\mathsurround=0pt\it $\cal G\mskip-2mu$eometry \&\ 
$\cal T\!\!$opology $\cal P\!$ublications}}  

\def\recd{{\small Received:\qua\receiveddate\ifx\reviseddate\relax
\else\qquad Revised:\qua\reviseddate\fi\par}} 

\def\lognumber#1{\def\thelognumber{#1}}
\def\volumenumber#1{\def\thevolumenumber{#1}}
\def\volumeyear#1{\def\thevolumeyear{#1}}
\def\papernumber#1{\def\thepapernumber{#1}}
\def\pagenumbers#1#2{\def\startpage{#1}\def\finishpage{#2}}
\def\published#1{\def\publishdate{#1}}

\def\received#1{\def\receiveddate{#1}}
\def\revised#1{\def\reviseddate{#1}}
\def\accepted#1{\def\accepteddate{#1}}

\let\\\par\let\thelognumber\relax\let\thevolumenumber\relax
\let\thepapernumber\relax\let\thevolumeyear\relax\let\startpage\relax
\let\finishpage\relax\let\publishdate\relax\let\receiveddate\relax
\let\reviseddate\relax\let\accepteddate\relax\let\theasciititle\relax
\let\theasciiauthors\relax
\let\theasciiabstract\relax

\let\theasciiemail\relax

\ifplaintex
\font\logobig=cmssbx10 scaled 3836
\font\logomed=cmssbx10 scaled 2557
\else
\font\logobig=cmssbx10 scaled 4200
\font\logomed=cmssbx10 scaled 2800
\fi

\long\def\makeagttitle{   
\count0=\startpage
\agt\hfill      
\hbox to 45truept{\vbox to 0pt{\vglue -13truept{\logomed A\kern -.37em{\logobig 
T}\kern -.38em G}\vss}\hss}
\break
{\small Volume \thevolumenumber\ (\thevolumeyear)
\startpage--\finishpage\nl
Published: \publishdate}

\vglue .25truein

{\parskip=0pt\leftskip 0pt plus
1fil\def\\{\par\smallskip}{\Large\bf\thetitle}\par\medskip} \vglue
0.05truein

{\parskip=0pt\leftskip 0pt plus 1fil\def\\{\par}{\sc\theauthors}
\par\medskip}%
 
\vglue 0.03truein 

{\small\leftskip 25truept\rightskip 25truept{\bf Abstract}\stdspace\theabstract

{\bf AMS Classification}\stdspace\theprimaryclass
\ifx\thesecondaryclass\relax\else; \thesecondaryclass\fi\par
{\bf Keywords}\stdspace \thekeywords\par}\vglue 7truept

}   

\ifplaintex
\font\phead=cmsl9 scaled 950
\font\pnum=cmbx10 scaled 913
\font\pfoot=cmsl9 scaled 950
\headline{\vbox to 0pt{\vskip -4.5mm\line{\small\phead\ifnum
\count0=\startpage ISSN 1472-2739 (on-line) 1472-2747 (printed)
\hfill {\pnum\folio}\else\ifodd\count0\def\\{ }%
\ifx\theshorttitle\relax\thetitle\else\theshorttitle\fi\hfill{\pnum\folio}
\else\def\\{ and }{\pnum\folio}\hfill\ifx\theshortauthors\relax\theauthors
\else\theshortauthors\fi\fi\fi}\vss}}
\footline{\vbox to 0pt{\vglue 0mm\line{\small\pfoot\ifnum\count0=\startpage
\copyright\ \gtp\hfill\else
\agt, Volume \thevolumenumber\ (\thevolumeyear)\hfill\fi}\vss}}
\else
\font=cmsl9 scaled 1050
\font\lnum=cmbx10 
\font=cmsl9 scaled 1050
\makeatletter
\def\@oddhead{{\small\lhead\ifnum\count0=\startpage ISSN 1472-2739 
(on-line) 1472-2747 (printed)\hfill {\lnum\number\count0}\else\ifodd\count0
\def\\{ }\ifx\theshorttitle\relax \thetitle \else\theshorttitle\fi\hfill
{\lnum\number\count0}\else\def\\{ and }{\lnum\number\count0}
\hfill\ifx\theshortauthors\relax 
\theauthors\else\theshortauthors\fi\fi\fi}}\def\@evenhead{\@oddhead}
\def\@oddfoot{\small\lfoot\ifnum\count0=\startpage\copyright\ \gtp\hfill\else
\agt, Volume \thevolumenumber\ (\thevolumeyear)\hfill\fi}
\def\@evenfoot{\@oddfoot}
\makeatother
\fi
\let\maketitlepage\makeagttitle
\let\makeshorttitle\maketitlepage
\let\maketitle\maketitlepage

\newwrite\gtoutfile
\long\gdef\makeheadfile{  
{\def\\{, }\def\s{ }
\immediate\openout\gtoutfile head.xxx
\immediate\write\gtoutfile{Proxy-for: \ifx\theasciiauthors\relax
\theauthors\else\theasciiauthors\fi\s<\ifx\theasciiemail\relax\theemail\else\theasciiemail\fi>}
\immediate\write\gtoutfile{\noexpand\\}
\immediate\write\gtoutfile{Authors: \ifx\theasciiauthors\relax
\theauthors\else\theasciiauthors\fi}
{\def\\{ }\immediate\write\gtoutfile{Title: \ifx\theasciititle\relax
\thetitle\else\theasciititle\fi}}
\immediate\write\gtoutfile{Subj-class: GT or SG, GR etc}
\immediate\write\gtoutfile{MSC-class: \theprimaryclass\ifx\thesecondaryclass\relax\else, \thesecondaryclass\fi}
\immediate\write\gtoutfile{Journal-ref: Algebr. Geom. Topol. \thevolumenumber\s
(\thevolumeyear) \startpage-\finishpage}
\immediate\write\gtoutfile{Comments: Published by Algebraic and
Geometric Topology at}
\immediate\write\gtoutfile{\s\s\s  http://www.maths.warwick.ac.uk/agt/AGTVol\thevolumenumber/agt-\thevolumenumber-\thepapernumber.abs.html}
\immediate\write\gtoutfile{\noexpand\\}
\immediate\write\gtoutfile{}
\ifx\theasciiabstract\relax
\immediate\write\gtoutfile{\theabstract}\else
\immediate\write\gtoutfile{\theasciiabstract}\fi
\immediate\write\gtoutfile{}
\immediate\write\gtoutfile{\noexpand\\}
\immediate\write\gtoutfile{}
\immediate\closeout\gtoutfile}}  

\def\maketitlepage{\makeagttitle\makeheadfile}
\let\makeshorttitle\maketitlepage
\let\maketitle\maketitlepage